
%

\input amstex
\loadbold
\documentstyle{amsppt}
\NoBlackBoxes

\pagewidth{32pc}
\pageheight{44pc}
\magnification=\magstep1

\def\ZZ{\Bbb{Z}}
\def\QQ{\Bbb{Q}}
\def\RR{\Bbb{R}}
\def\CC{\Bbb{C}}

\def\OO{\Cal{O}}
\def\Spec{\operatorname{Spec}}
\def\Pic{\operatorname{Pic}}
\def\zeros{\operatorname{div}}
\def\GL{\operatorname{GL}}
\def\rank{\operatorname{rk}}
\def\rest#1#2{\left.{#1}\right\vert_{{#2}}}
\def\achern#1#2{\widehat{c}_{#1}(#2)}
\def\achch#1#2{\widehat{\operatorname{ch}}_{#1}(#2)}
\def\achow{\widehat{\operatorname{CH}}}
\def\adeg{\widehat{\operatorname{deg}}}
\def\chch#1#2{\operatorname{ch}_{#1}(#2)}
\def\SStable{\overline{\boldkey M}}

\topmatter
\title
Hodge index theorem for \\
arithmetic cycles of codimension one
\endtitle
\rightheadtext{}
\author Atsushi Moriwaki \endauthor
\leftheadtext{}
\address
Department of Mathematics, Faculty of Science,
Kyoto University, Kyoto, 606-01, Japan
\endaddress
\curraddr
Department of Mathematics, University of California,
Los Angeles, 405 Hilgard Avenue, Los Angeles, California 90024, USA
\endcurraddr
\email moriwaki\@math.ucla.edu \endemail
\date March, 1994 \enddate
\thanks
This is a revised version of my previous paper
``Hodge index theorem on arithmetic varieties''.
\endthanks
\endtopmatter

\document

\subhead 0. Introduction
\endsubhead

Let $f : X \to \Spec(\ZZ)$ be a $(d+1)$-dimensional
regular arithmetic variety over
$\Spec(\ZZ)$, i.e. $X$ is regular, $X$ is projective and
flat over $\Spec(\ZZ)$ and $d = \dim f$. Let $H$ be an $f$-ample line bundle on
$X$
and $k$ a Hermitian metric of $H$.
Here we consider a homomorphism
$$
 L : \achow^p(X)_{\RR} \to \achow^{p+1}(X)_{\RR}
$$
defined by $L(x) = x \cdot \achern{1}{H, k}$.
In \cite{GS}, H. Gillet and C. Soul\'{e} conjectured that

\proclaim{Arithmetic Analogues of Grothendieck's Standard Conjectures}
For a suitable choice of $k$, if $2p \leq d+1$, then
\roster
\item "(a)" The homomorphism
$L^{d+1-2p} : \achow^p(X)_{\RR} \to \achow^{d+1-p}(X)_{\RR}$
is bijective, and

\item "(b)" If $x \in \achow^p(X)_{\RR}$, $x \not= 0$ and
$L^{d+2-2p}(x) = 0$, then
$(-1)^p \adeg(x \cdot L^{d+1-2p}(x)) > 0$.
\endroster
\endproclaim

\noindent
For example, K. K\"{u}nnemann \cite{Ku} proved that if $X$ is
a projective space, then the conjecture is true.
In this note, we would like to prove the following partial answer
of the above conjecture for general arithmetic varieties.

\proclaim{Theorem A}
Assume that $d \geq 1$ and $(H, k)$
is arithmetically ample, i.e. $(H, k)$ satisfies the following
conditions:
\roster
\item "(i)" $H$ is $f$-ample.

\item "(ii)" The Chern form $c_1(H_{\infty}, k_{\infty})$ gives a K\"{a}hler
form
            on the infinite fiber $X_{\infty}$.

\item "(iii)" For every irreducible horizontal subvariety $Y$
            (i.e. $Y$ is flat over $\Spec(\ZZ)$), the height
            $\achern{1}{\rest{(H, k)}{Y}}^{\dim Y}$ of $Y$ is positive.
\endroster
Then we have the following:
\roster
\item "(1)" $L^{d-1} : \achow^1(X)_{\RR} \to \achow^{d}(X)_{\RR}$
is injective.

\item "(2)" If $x \in \achow^1(X)_{\RR}$, $x \not= 0$ and
$L^{d}(x) = 0$, then $\adeg(x L^{d-1}(x)) < 0$.
\endroster
\endproclaim

\noindent
Theorem~A is a consequence of the following higher dimensional
generalization of Faltings-Hriljac's Hodge index theorem
on arithmetic surfaces (cf. \cite{Fa} and \cite{Hr}).

\proclaim{Theorem B}
Assume that $d \geq 1$ and $(H, k)$ is arithmetically ample. Let
$X \overset f' \to \longrightarrow \Spec(O_K) \to \Spec(\ZZ)$
be the Stein factorization of $f : X \to \Spec(\ZZ)$
and $X_K$ the generic fiber of $f'$,
where $O_K$ is the ring of integers of an algebraic number field $K$.
Let $z : \achow^1(X) \to \operatorname{CH}^1(X)$ be the
canonical homomorphism defined by $z(D, g) = D$.
If $x \in \achow^1(X)$ and
$\left(\rest{z(x)}{X_K} \cdot \left(\rest{H}{X_K}\right)^{d-1}\right) = 0$,
then
$$
\adeg(x^2 \cdot \achern{1}{H, k}^{d-1}) \leq 0.
$$
Moreover, equality holds if and only if
there are a positive integer $n$ and $y \in \achow^1(\Spec(O_K))$
such that $nx = f^*(y)$.
\endproclaim

\subhead 1. Proof of Theorem~B
\endsubhead

In this section, we would like to give the proof of Theorem~B.
An advantage to use arithmetical ampleness of the
Hermitian line bundle $(H, k)$ due to S. Zhang \cite{Zh}
is that a higher multiple of it produces a lot of good sections
(cf. \cite{Zh} and \cite{Mo2}), so that we can proceed induction
on $d = \dim f$. However, regularity of $X$ doesn't preserve
by induction step in general. Here we consider the following
weaker version on general arithmetic varieties.

\proclaim{Theorem 1.1}
Let $K$ be an algebraic number field and $O_K$ the ring of integers of $K$.
Let $f : X \to \Spec(O_K)$ be an arithmetic variety such that
$d = \dim f \geq 1$ and $X_K$ is smooth and geometrically irreducible.
Let $(H, k)$ be an arithmetically ample Hermitian line bundle on $X$, i.e.
$(H, k)$ is a Hermitian line bundle with the following properties:
\roster
\item "(i)" $H$ is $f$-ample.

\item "(ii)" The Chern form $c_1(H_{\sigma}, k_{\sigma})$ gives
            a K\"{a}hler form on $X_{\sigma}$ for all $\sigma \in K(\CC)$.

\item "(iii)" For every irreducible horizontal subvariety $Y$
            (i.e. $Y$ is flat over $\Spec(O_K)$), the height
            $\achern{1}{\rest{(H, k)}{Y}}^{\dim Y}$ of $Y$ is positive.
\endroster
Let $D$ be a Cartier divisor on $X$ and $g_{\sigma}$
a Green current of $D_{\sigma}$ on each $\sigma \in K(\CC)$.
If $(D_K \cdot H_{K}^{d-1}) = 0$,
then
$$
\adeg\left(\left(D, \sum g_{\sigma}\right)^2
\cdot \achern{1}{H, k}^{d-1}\right) \leq 0.
$$
Moreover, if equality holds, then
there is a positive integer $n$,
a Cartier divisor $Z$ on $X$ and constants
$\{ g'_{\sigma} \}_{\sigma \in K(\CC)}$ such that
the support of $Z$ is vertical and the class of
$(Z, \sum g'_{\sigma})$ is equal to
the class of $n(D, \sum g_{\sigma})$ in $\achow^1(X)$.
In particular, if equality holds, then
$\OO_{X_K}(D_K)$ is a torsion of $\Pic(X_K)$.
\endproclaim

\demo{Proof}
First of all, we prepare two lemmas.

\proclaim{Lemma 1.1.1}
Let $X$ be a $d$-dimensional compact K\"{a}hler manifold
with a K\"{a}hler form $\Phi$ and
$\varphi$ a real valued smooth function on $X$.
Then,
$$
    \int_X \varphi dd^c(\varphi) \Phi^{d-1} \leq 0.
$$
Moreover, equality holds if and only if $\varphi$ is a constant.
\endproclaim

\demo{Proof}
Since ${\displaystyle dd^c = \frac{\sqrt{-1}}{2\pi}\partial\bar{\partial}}$
and $d(\varphi\bar{\partial}(\varphi))
= \partial(\varphi)\bar{\partial}(\varphi) +
\varphi \partial\bar{\partial}(\varphi)$, by Stokes' theorem,
we have
$$
\int_X \varphi dd^c(\varphi) \Phi^{d-1} =
- \frac{\sqrt{-1}}{2\pi} \int_X \partial(\varphi)\bar{\partial}(\varphi).
$$
Here let $\theta^1, \ldots, \theta^d$ be a local unitary frame of
$\Omega^1_X$ with
$\Phi = \sqrt{-1} \sum_{i} \theta^{i} \wedge \bar{\theta}^{i}$.
We set $\partial(\varphi) =
\sum_{i} a_{i} \theta^{i}$. Then,
$\bar{\partial}(\varphi) = \overline{\partial(\varphi)}
= \sum_{i} \bar{a}_{i} \bar{\theta}^{i}$.
Therefore,
$$
- \frac{\sqrt{-1}}{2\pi} \partial(\varphi)\bar{\partial}(\varphi)
= \frac{-1}{2 \pi} \sum_{i=1}^d |a_i|^2 \Phi^d.
$$
Thus, we have
$$
    \int_X \varphi dd^c(\varphi) \Phi^{d-1} \leq 0.
$$
Moreover, equality hold if and only if $\partial(\varphi) = 0$.
Here, since $\varphi$ is real valued, $\partial(\varphi) = 0$
implies that $\varphi$ is a constant.
\qed
\enddemo

\proclaim{Lemma 1.1.2}
Let $X$ be a $d$-dimensional K\"{a}hler manifold
with a K\"{a}hler form $\Phi$ and
$\omega$ a smooth $(1,1)$-form on $X$ such that $\bar{\omega} = - \omega$
and $\omega \wedge \Phi^{d-1} = 0$.
Then, there is a real valued smooth function $u$ on $X$ with
the following properties.
\roster
\item "(1)" $\omega^2 \wedge \Phi^{d-2} = u \Phi^{d}$.
\item "(2)" $u(x) \leq 0$ for all $x \in X$.
\item "(3)" $u(x) = 0$ for all $x \in X$ if and only if $\omega = 0$.
\endroster
\endproclaim

\demo{Proof}
Let $\theta^1, \ldots, \theta^d$ be a local unitary frame of
$\Omega^1_X$ with
$\Phi = \sqrt{-1} \sum_{i} \theta^{i} \wedge \bar{\theta}^{i}$.
We set $\omega = \sum_{i, j} a_{ij} \theta^{i} \wedge \bar{\theta}^{j}$.
Then, $\bar{\omega} = - \omega$ implies that $a_{ji} = -\bar{a}_{ij}$.
Moreover, since
$$
\omega \wedge \Phi^{d-1} = - \sqrt{-1} (a_{11} + \cdots + a_{dd}) \Phi^d,
$$
we have $a_{11} + \cdots + a_{dd} = 0$.
On the other hand, by an easy calculation,
$$
d(d-1) \omega^2 \wedge \Phi^{d-2} =
\left(\sum_{i, j} a_{ij}a_{ji} - a_{ii}a_{jj}
\right) \Phi^d.
$$
Therefore, we get
$$
\omega^2 \wedge \Phi^{d-2} =
\frac{-1}{d(d-1)} \left(\sum_{i, j} |a_{ij}|^2 \right) \Phi^{d}.
$$
Hence, if we set
${\displaystyle u = \frac{-1}{d(d-1)} \sum_{i, j} |a_{ij}|^2}$,
the lemma is obtained because $\sum_{i, j} |a_{ij}|^2$ is
independent of the choice of $\theta^1, \ldots, \theta^d$.
\qed
\enddemo

Let us start of the proof of Theorem~1.1.
We will prove it by induction on $d$.
First, we consider the case $d=1$.
In this case, taking a desingularization of
$X$, we may assume that $X$ is regular.
Thus, our theorem can be derived from
Faltings-Hriljac's Hodge index theorem (cf. \cite{Fa} and \cite{Hr}).

\medskip
Next, we assume $d \geq 2$.
We set $x = (D, \sum g_{\sigma})$ and
$L = \OO_{X}(D)$. Let $h_{\sigma}$ be an Einstein-Hermitian
metric of $L_{\sigma}$ with respect to $c_1(H_{\sigma}, h_{\sigma})$.
Let $s$ be a rational section of $L$ with
$\operatorname{div}(s) = D$.
Here we consider an arithmetic cycle
$$
y = \left(
D, \sum_{\sigma \in K(\CC)} -\log(h_{\sigma}(s_{\sigma}, s_{\sigma}))
\right).
$$
Since $g_{\sigma}$ and
$-\log(h_{\sigma}(s_{\sigma}, s_{\sigma}))$ are
Green currents of the same $D_{\sigma}$,
there is a real valued smooth function
$\phi_{\sigma}$ on each $X_{\sigma}$ such that
$x = y + a(\sum_{\sigma \in K(\CC)} \phi_{\sigma})$ in $\achow^1(X)$.
Then,
it is easy to see that
$$
\adeg(x^2 \cdot \achern{1}{H, k}^{d-1}) =
\adeg(y^2 \cdot \achern{1}{H, k}^{d-1}) + \frac{1}{2} \sum_{\sigma \in K(\CC)}
\int_{X_{\sigma}} \phi_{\sigma} dd^c(\phi_{\sigma}) c_1(H_{\sigma},
k_{\sigma})^{d-1}
$$
because
$c_1(L_{\sigma}, h_{\sigma}) c_1(H_{\sigma}, k_{\sigma})^{d-1} = 0$.
Therefore, by Lemma~1.1.1,
$$
\adeg(x^2 \cdot \achern{1}{H, k}^{d-1}) \leq
\adeg(y^2 \cdot \achern{1}{H, k}^{d-1})
$$
and equality holds if and only if $\phi_{\sigma}$ is a constant
for each $\sigma \in K(\CC)$.
On the other hand, by virtue of \cite{Mo2, Theorem~4.2 and Theorem~5.2},
for a sufficiently large $m$, there is a section $t \in H^0(X, H^m)$
with the following properties:
\roster
\item "i)" $\zeros(t)_K$ is smooth and geometrically irreducible.

\item "ii)" If $\zeros(t) = Y + a_1 F_1 + \cdots + a_s F_s$
is the irreducible decomposition such that
$Y$ is horizontal and $F_i$'s are vertical, then
$F_i$'s are smooth fibers.

\item "iii)" $D$ and $\zeros(t)$ has no common irreducible component.

\item "iv)" $\sup_{x \in X_{\sigma}}
\left( ||t_{\sigma}||_{k_{\sigma}^m}(x) \right) < 1$ for all $\sigma \in
K(\CC)$.
\endroster
(Note that $H^1(X_K, H_K^{-m}) = 0$ guarantees
geometrical irreducibility of $\zeros(t)_K$.)
Since $(\rest{D}{F_i}^2 \cdot \rest{H}{F_i}^{d-2}) \leq 0$
by the geometric Hodge index theorem, we obtain
$$
\align
\adeg(y^2 \cdot \achern{1}{H^m, k^m}^{d-1})  & =
      \adeg(\rest{y}{Y}^2 \cdot \achern{1}{\rest{(H^m, k^m)}{Y}}^{d-2})
+ \sum a_i m (\rest{D}{F_i}^2 \cdot \rest{H}{F_i}^{d-2}) \\
& \qquad
 - \sum_{\sigma \in K(\CC)}
   \int_{X_{\sigma}} \log(||t_{\sigma}||_{k_{\sigma}^m})
    c_1(L_{\sigma}, h_{\sigma})^2 c_1(H^m_{\sigma}, k^m_{\sigma})^{d-2} \\
& \leq
      \adeg(\rest{y}{Y}^2 \cdot \achern{1}{\rest{(H^m, k^m)}{Y}}^{d-2}) \\
& \qquad
 - \sum_{\sigma \in K(\CC)}
   \int_{X_{\sigma}} \log(||t_{\sigma}||_{k_{\sigma}^m})
    c_1(L_{\sigma}, h_{\sigma})^2 c_1(H^m_{\sigma}, k^m_{\sigma})^{d-2}.
\endalign
$$
Since $(L_{\sigma}, h_{\sigma})$ is Einstein-Hermitian, by Lemma~1.1.2,
there is a real-valued smooth function $u_{\sigma}$ on $X_{\sigma}$
with the following properties:
\roster
\item "(1)" $c_1(L_{\sigma}, h_{\sigma})^2 c_1(H_{\sigma}, k_{\sigma})^{d-1}
= u_{\sigma} c_1(H_{\sigma}, k_{\sigma})^{d}$.
\item "(2)" $u_{\sigma}(x) \leq 0$ for all $x \in X_{\sigma}$.
\item "(3)" $u_{\sigma}(x) = 0$ for all $x \in X_{\sigma}$ if and only if
$(L_{\sigma}, h_{\sigma})$ is flat.
\endroster
Therefore, we have
$$
\adeg(y^2 \cdot \achern{1}{H, k}^{d-1}) \leq
      \adeg(\left(\rest{y}{Y}\right)^2 \cdot \achern{1}{\rest{(H,
k)}{Y}}^{d-2}).
$$
Hence, by hypothesis of induction, we get our inequality.

\medskip
Finally, we consider equality condition.
We assume $\adeg(x^2 \cdot \achern{1}{H, k}^{d-1}) = 0$.
Then, if we trace back the above proof carefully, we can see
\roster
\item "(a)" $\phi_{\sigma}$ is a constant for each $\sigma \in K(\CC)$.

\item "(b)" $(L_{\sigma}, h_{\sigma})$ is flat for each $\sigma \in K(\CC)$.

\item "(c)" $\rest{L_K}{Y_K}$ is a torsion of $\Pic(Y_K)$.
\endroster
By (b), $L_{\CC}$ is given by
a representation $\rho : \pi_1(X_{\CC}) \to \CC^{*}$
of the fundamental group of $X_{\CC}$.
(c) implies that the image of
$\pi_1(Y_{\CC}) \to \pi_1(X_{\CC}) \to \CC^{*}$ is finite.
On the other hand, by Lefschetz theorem (cf. Theorem 7.4 in \cite{Mi}),
$\pi_1(Y_{\CC}) \to \pi_1(X_{\CC})$ is surjective.
Thus, the image of $\rho : \pi_1(X_{\CC}) \to \CC^{*}$ is also finite.
Therefore, there is a positive integer $n$ with
$L_{\CC}^n \simeq\OO_{X_{\CC}}$. Thus,
$$
\dim_K H^0(X_K, L_K^n) = \dim_{\CC} H^0(X_K, L_K^n) \otimes \CC =
\dim_{\CC} H^0(X_{\CC}, L_{\CC}^n) = 1.
$$
Hence, since $(L_K \cdot H_K^{d-1}) = 0$, we have $L_K^n \simeq \OO_{X_K}$.
Thus, there is a rational section $s'$ of $L^n$
with $s'_K = 1$.
We set $Z= \operatorname{div}(s')$ and
$g'_{\sigma} =  -\log(h_{\sigma}^n(s', s')) + n \phi_{\sigma}$.
Then, the support of $Z$ is vertical.
Moreover, since $h_{\sigma}^n$ is a flat metric of $\OO_{X_{\sigma}}$,
$h_{\sigma}^n(s', s')$ must be a constant.
Therefore, $(Z, \sum g'_{\sigma})$ is our desired cycle.
\qed
\enddemo

\bigskip
\demo{Proof of Theorem B}
Since $f'_* \OO_X = O_K$, $X_K$ is geometrically irreducible.
So the inequality is an immidiate consequence of Theorem~1.1.

We need to consider the precise equality condition.
Clearly, if there are a positive integer $n$ and
$y \in \achow^1(\Spec(O_K))$ such that $nx = {f'}^*(y)$, then
$\adeg(x^2 \cdot \achern{1}{H, k}^{d-2}) = 0$.
Conversely we assume $\adeg(x^2 \cdot \achern{1}{H, k}^{d-2}) = 0$.
Then, by Theorem~1.1,
there are a positive integer $n_1$ and an arithmetic cycle
$(Z, \sum_{\sigma \in K(\CC)} g_{\sigma})$
such that $Z$ is vertical with respect to $f'$,
$g_{\sigma}$'s are constant and
$n_1x$ is equal to the class of
$(Z, \sum_{\sigma \in K(\CC)} g_{\sigma})$
in $\achow^1(X)$. Then,
$$
\adeg((n_1x)^2 \cdot \achern{1}{H, k}^{d-2}) =
(Z^2 \cdot H^{d-1}) = 0.
$$
Here, we need the following lemma.

\proclaim{Lemma 1.3}
Let $X$ be a regular scheme, $R$ a discrete valuation ring,
$f : X \to \Spec(R)$ a projective morphism with $f_* \OO_X = R$, and
$H$ an $f$-ample line bundle on $X$.
Let $X_o$ be the central fiber of $f$ and
$(X_o)_{\operatorname{red}} = X_1 + \cdots + X_n$ the irreducible decomposition
of $(X_o)_{\operatorname{red}}$.
We consider a vector space $V = \bigoplus_{i=1}^n \QQ X_i$
generated by $X_i$'s and the natural pairing $(\ , \ ) : V \times V \to \QQ$
defined by
$$
    (D_1, D_2) = (D_1 \cdot D_2 \cdot H^{d-1}),
$$
where $d = \dim f$ and $\cdot$ is the intersection product.
Then, we have $(D, D) \leq 0$ for all $D \in V$ and
equality holds if and only if $D \in \QQ X_o$.
\endproclaim

\demo{Proof}
For example, see (i)' of Lemma (2.10) in Chap. I of \cite{BPV}.
\qed
\enddemo

\noindent
By the above lemma, there is a positive integer $n_2$
and a cycle $T$ on $\Spec(\OO_K)$ such that
$n_2 Z = {f'}^*(T)$.
Therefore, if we set $y = (T, \sum_{\sigma \in K(\CC)} n_2 g_{\sigma})$,
then $n_1 n_2 x = {f'}^*(y)$.
\qed
\enddemo

\subhead 2. Proof of Theorem A
\endsubhead

Let us begin the proof of Theorem~A, This is
an easy corollary of Theorem~B.

\bigskip
(1) Let us see that (2) implies (1).
Assume that $L^{d-1}(x) = 0$. Then, $L^d(x) = 0$.
Thus if $x \not= 0$, then $\adeg(x L^{d-1}(x)) < 0$ by (2).
This is a contradiction. Therefore, $x = 0$.

\bigskip
(2) Let $X \overset f' \to \longrightarrow \Spec(O_K) \to \Spec(\ZZ)$
be the Stein factorization of $f : X \to \Spec(\ZZ)$.
In the following arguments, the subscript $K$
means the restriction to the generic fiber of $f'$.

Since $x$ can be approximated by points $y \in \achow^1(X)_{\QQ}$
with $L^d(y) = 0$, we may assume that $x \in \achow^1(X)_{\QQ}$.
Let $t$ be a rational number with
$(z(x)_K + t H_K \cdot H_K^{d-1}) = 0$.
Replacing $x$ by $mx$, we may assume that
$x \in \achow^1(x)$ and $t \in \ZZ$.
We set $y = x + t \achern{1}{H, k}$.
Then, $(z(y)_K \cdot H_K^{d-1}) = 0$. Thus, by Theorem~B,
we have $\adeg(y^2 \cdot \achern{1}{H, k}^{d-1}) \leq 0$.
Therefore, since $L^d(x) = 0$, we get
$$
\adeg(x^2 \cdot \achern{1}{H, k}^{d-1}) + (t)^2
\adeg(\achern{1}{H, k}^{d+1}) \leq 0.
$$
Hence, $\adeg(x^2 \cdot \achern{1}{H, k}^{d-1}) \leq 0$.
Here, we assume that $\adeg(x^2 \cdot \achern{1}{H, k}^{d-1}) = 0$.
Then, $t = 0$. Thus, $(z(x)_K \cdot H_K^{d-1}) = 0$.
So, by Theorem~B,
there is a positive integer $n$ and
$u \in \achow^1(\Spec(O_K))$ such that $nx = {f'}^*(u)$.
We know $nx \cdot \achern{1}{H, k}^{d} = 0$, which implies
$u \cdot f'_*(\achern{1}{H, k}^{d}) = 0$. Therefore,
$u = 0$ in $\achow^1(\Spec(O_K))_{\QQ}$
because $f'_*(\achern{1}{H, k}^{d}) = (H_K^d)[\Spec(O_K)]$.
Thus, $x = 0$ in $\achow^1(X)_{\QQ}$.
This is a contradiction.
Hence, we get $\adeg(x^2 \cdot \achern{1}{H, k}^{d-1}) < 0$.

\subhead 3. Variants of Theorem B
\endsubhead

In this section, we will study variants of Theorem~B or Theorem~1.1.
The following theorem is a generalization of Theorem~1.1 to a higher
rank vector bundle.

\proclaim{Theorem 3.1}
Let $K$ be an algebraic number field and $O_K$ the ring of integers.
Let $f : X \to \Spec(O_K)$ be an arithmetic variety
and $(H, k)$ an arithmetically ample Hermitian line bundle on $X$.
Assume that $d = \dim f \geq 1$ and $X_K$ is smooth and
geometrically irreducible.
Let $(E, h)$ be a Hermitian vector bundle on $X$ such that
$E_{\overline{\QQ}}$ is semi-stable with respect to $H_{\overline{\QQ}}$
and $(c_1(E_K) \cdot c_1(H_K)^{d-1}) = 0$. Then, we have
$$
\achch{2}{E, h} \cdot \achern{1}{H, k}^{d-1} \leq 0.
$$
Moreover, if the equality holds, then $h_{\sigma}$ is Einstein-Hermitian
with respect to a K\"{a}hler form
$\Omega_{\sigma} = c_1(H_{\sigma}, k_{\sigma})$
and $E_{\sigma}$ is flat for every $\sigma \in K(\CC)$.
\endproclaim

\demo{Proof}
Let $r$ be the rank of $E$. Since
$$
\achch{2}{E, h} = \frac{1}{2}\achern{1}{E, h}^2 - \achern{2}{E, h},
$$
we have
$$
\align
\achch{2}{E, h} \cdot \achern{1}{H, k}^{d-1} & =
\frac{1}{2r} \achern{1}{E, h}^2 \cdot \achern{1}{H, k}^{d-1} \\
&\phantom{=}\quad -
\left\{ \achern{2}{E, h} - \frac{r-1}{2r} \achern{1}{E, h}^2 \right\} \cdot
\achern{1}{H, k}^{d-1}.
\endalign
$$
By Lemma~8.2 of \cite{Mo1}, $E_{\sigma}$ is semistable with respect
to $H_{\sigma}$. Thus the main theorem in \cite{Mo2} implies that
$$
\left\{ \achern{2}{E, h} - \frac{r-1}{2r} \achern{1}{E, h}^2 \right\} \cdot
\achern{1}{H, k}^{d-1} \geq 0.
$$
On the other hand, by Theorem~1.1,
$\achern{1}{E, h}^2 \cdot \achern{1}{H, k}^{d-1} \leq 0$.
Therefore, we have $\achch{2}{E, h} \cdot \achern{1}{H, k}^{d-1} \leq 0$.

Next we consider equality condition.
We assume that $\achch{2}{E, h} \cdot \achern{1}{H, k}^{d-1} = 0$.
First of all, by equality condition of the main theorem of \cite{Mo2},
$E_{\sigma}$ is flat for every $\sigma \in K(\CC)$.
Let $h'$ be an Einstein-Hermitian metric of $E$.
Then, by Lemma 6.1 of \cite{Mo1},
$$
(\achch{2}{E, h} - \achch{2}{E, h'}) \cdot \achern{1}{H, k}^{d-1}
= -\frac{(d-1)!}{4 \pi} \sum_{\sigma \in K(\CC)}
DL(E_{\sigma}, h_{\sigma}, h'_{\sigma}),
$$
where $DL$ is the Donaldson's Lagrangian.
Therefore, we have
$$
\sum_{\sigma \in K(\CC)} DL(E_{\sigma}, h_{\sigma}, h'_{\sigma}) \leq 0.
$$
On the other hand, since $h'$ is Einstein-Hermitian,
we get $DL(E_{\sigma}, h_{\sigma}, h'_{\sigma}) \geq 0$
for all $\sigma \in K(\CC)$.
Hence $DL(E_{\sigma}, h_{\sigma}, h'_{\sigma}) = 0$ for all
$\sigma \in K(\CC)$. Thus $h_{\sigma}$ is Einstein-Hermitian for
all $\sigma \in K(\CC)$.
\qed
\enddemo

In the case where $\rank E = 1$, Theorem 1.1 says that
if $\achch{2}{E, h} \cdot \achern{1}{H, k}^{d-1} = 0$, then
$E_K$ is a torsion element of $\Pic^0(X_K)$.
So we might expect a stronger property of $(E, h)$ than flatness.
Here we introduce one notation. Let $M$ be a complex manifold and
$F$ a flat vector bundle of rank $r$ on $M$.
Let $\rho_F : \pi_1(M) \to \GL_r(\CC)$ be the representation of the
fundamental group of $M$ arising from the flat vector bundle $F$.
$F$ is said to be {\it of torsion type} if the image of $\rho_F$ is finite.

\proclaim{Proposition 3.2}
Let $K$ be an algebraic number field and $O_K$ the ring of integers.
Let $f : X \to \Spec(O_K)$ be an arithmetic variety,
$H$ an $f$-ample line bundle on $X$ and
$k$ a Hermitian metric of $H$.
Assume that $d = \dim f \geq 1$ and $X_K$ is smooth and
geometrically irreducible.
Let $(E, h)$ be a Hermitian vector bundle of rank $r$ on $X$ such that
$(E_{\sigma}, h_{\sigma})$ is flat for each $\sigma \in K(\CC)$
and $\achch{2}{E, h} \cdot \achern{1}{H, k}^{d-1} = 0$.
Let $\rho_{E_{\CC}} : \pi_1(X_{\CC}) \to \GL_r(\CC)$ be the representation
of the fundamental group of $X_{\CC}$
arising from the flat vector bundle $E_{\CC}$.
If the image of $\rho_{E_{\CC}}$ is abelian,
then $E_{\sigma}$ is of torsion type
for all $\sigma \in K(\CC)$.
\endproclaim

\demo{Proof}
We prove it by induction on $\dim X$.
First, we consider the case $d = 1$.
Since the representation $\rho_{E_\CC}$ is abelian,
we have the decomposition
$\rho_{E_\CC} = \rho_1 \oplus \cdots \oplus \rho_r$
such that $\dim \rho_i = 1$ for all $i$.
Therefore, there are flat line bundles
$L'_1, \ldots, L'_r$ on $X_{\CC}$ such that
$E_{\CC} = L'_1 \oplus \cdots \oplus L'_r$.
Thus, by an easy descent, we can find line bundles
$L_1, \ldots, L_r$ on $X_{\overline{\QQ}}$ such that
$E_{\overline{\QQ}} = L_1 \oplus \cdots \oplus L_r$ and $\deg(L_i) = 0$
for all $i$. Thus, by Proposition~10.8 in \cite{Mo1},
we have our assertion.

Next, we may assume that $d \geq 2$.
Replacing $H$ by a higher multiple $H^m$ of $H$,
we may assume that there is a section $\phi \in H^0(X, H)$
with the following properties:
\roster
\item "i)" $\zeros(\phi)_K$ is smooth and geometrically irreducible.

\item "ii)" If $\zeros(\phi) = Y + a_1 F_1 + \cdots + a_s F_s$
is the irreducible decomposition such that
$Y$ is horizontal and $F_i$'s are vertical,
then $F_i$'s are smooth fibers.
\endroster
Since $(E_{\sigma}, h_{\sigma})$ is flat for each
$\sigma \in K(\CC)$,
we have $(\chch{2}{E} \cdot F_i \cdot H^{d-2}) = 0$ and
$\chch{2}{E_{\sigma}, h_{\sigma}}$ is zero as differential form for
every $\sigma \in K(\CC)$.
Thus we have
$$
\achch{2}{E, h} \cdot \achern{1}{H, k}^{d-1}  =
      \achch{2}{\rest{(E, h)}{Y}} \cdot \achern{1}{\rest{(H, k)}{Y}}^{d-2}.
$$
Let $\rho_{\rest{E_{\CC}}{Y_{\CC}}} : \pi_1(Y_{\CC}) \to \GL_r(\CC)$
be the representation arising from $\rest{E_{\CC}}{Y_{\CC}}$.
Since $\rho_{\rest{E_{\CC}}{Y_{\CC}}}$ is the composition
of $\pi_1(Y_{\CC}) \to \pi(X_{\CC})$ and
$\rho_{E_{\CC}} : \pi_1(X_{\CC}) \to \GL_r(\CC)$,
the image of $\rho_{\rest{E_{\CC}}{Y_{\CC}}}$ is also abelian.
Thus, by hypothesis of induction, $\rest{E_{\sigma}}{Y_{\sigma}}$
is of torsion type for every $\sigma \in K(\CC)$.
On the other hand, by Lefschetz theorem,
$\pi_1(Y_{\sigma}) \to \pi_1(X_{\sigma})$ is surjective. Hence,
$E_{\sigma}$ is also of torsion type for every $\sigma \in K(\CC)$.
\qed
\enddemo

Finally, we will pose two questions.
Let $f : X \to \Spec(O_K)$ be a $(d+1)$-dimensional arithmetic
variety, $(H, k)$ an arithmetically ample Hermitian line bundle on $X$,
and $(E, h)$ a Hermitian vector bundle on $X$ such that
$E_{\overline{\QQ}}$ is semistable with respect to $H_{\overline{\QQ}}$
and $(c_1(E_K) \cdot c_1(H_K)^{d-1}) = 0$.
An interesting problem is to find stronger equality conditions for
$$
\achch{2}{E, h} \cdot \achern{1}{H, k}^{d-1} \leq 0.
$$
Theorem~3.1 says that if
$\achch{2}{E, h} \cdot \achern{1}{H, k}^{d-1} = 0$, then
at least $E_{\sigma}$ is flat for every $\sigma \in K(\CC)$.
Optimistically, one may pose the following question:

\definition{Question 3.3}
If $\achch{2}{E, h} \cdot \achern{1}{H, k}^{d-1} = 0$,
is $E_{\sigma}$ of torsion type for every $\sigma \in K(\CC)$ ?
\enddefinition

\noindent
By Proposition~3.2, if $\pi_1(X_{\CC})$ is abelian or $\rank E = 1$, we have
an affirmative answer of the above question.
Moreover, if we carefully trace back the proof in Proposition~3.2,
Question~3.3 can be reduced to the case $d = 1$.
So from now on, we assume that $d = 1$.
Let $\SStable_{X_K/K}(r, 0)$ be
the moduli scheme of semistable vector bundles
on $X_K$ with rank $r$ and degree $0$. Let $h$ be a height function on
$\SStable_{X_K/K}(r, 0)$ arising from some ample
line bundle on $\SStable_{X_K/K}(r, 0)$. Our next question is

\definition{Question 3.4}
Are there constants $A$ and $B$ with the following properties ?
\roster
\item "(1)" $A, B \in \RR$ and $A > 0$.

\item "(2)" For all semistable Hermitian vector bundle $(E, h)$ on $X$
with rank $r$ and degree $0$, we have
$$
         h(E_K) \leq \frac{-A}{[K : \QQ]} \achch{2}{E, h} + B
$$
\endroster
\enddefinition

\noindent
In some sense, Question~3.4 is related to Question~3.3.
For, if $\achch{2}{E, h} = 0$ and Question~3.4 holds,
then the height of $E_K$ is bounded.
So $E_K$ should have some simple structure.


\enddocument